\documentclass[10pt,twoside]{article}
\usepackage{lineno}
\usepackage{CJK}
\usepackage{graphicx}
\usepackage{amsfonts}
\usepackage{fancyhdr}
\usepackage{titlesec}
\usepackage{cite}
\usepackage{ifthen}
\usepackage{pifont}
\usepackage{stmaryrd}
\usepackage{setspace}
\usepackage{indentfirst}
\usepackage{amsmath,amssymb,amscd,bbm,amsthm,mathrsfs,dsfont}
\usepackage{bm}
\usepackage{authblk}
\input amssym.def

\newtheorem{conjecture}{Conjecture}[section]
\newtheorem{theorem}{Theorem}[section]

\newtheorem{lemma}{Lemma}[section]
\newtheorem{remark}{Remark}[section]

\newtheorem*{pf}{Proof}
\theoremstyle{plain}

\newboolean{first}
\setboolean{first}{true}
\begin{document}
%\linenumbers  % ÆôÓÃÐкÅ
\title{\bf \Large Injective edge-coloring of claw-free graphs with maximum degree 4
\footnotetext[1]{
\textbf{The first author's research is supported by NSFC (No. 12171436)}
}}

\author{Danjun Huang$^*$}
\author{Yuqian Guo}
\affil{\normalsize{Department of Mathematics, Zhejiang Normal University, Jinhua
321004, China}}
\date{}
\maketitle

\noindent
\textbf{Abstract}\quad An injective $k$-edge-coloring of a graph $G$ is a mapping
$\phi$: $E(G)\rightarrow\{1,2,...,k\}$, such that $\phi(e)\ne\phi(e')$ if edges $e$
and $e'$ are at distance two, or are in a triangle. The smallest integer $k$ such that
$G$ has an injective $k$-edge-coloring is called the injective chromatic index of $G$,
denoted by $\chi_i'(G)$. A graph is called claw-free if it has no induced subgraph
isomorphic to the complete bipartite graph $K_{1,3}$. In this paper, we show that
$\chi_i'(G)\le 13$ for every claw-free graph $G$ with $\Delta(G)\leq 4$, where
$\Delta(G)$ is the maximum degree of $G$.
\bigskip

\bigskip
\bigskip
\section{Introduction}
\par  Only simple and finite graphs are considered in this paper. We use $V(G)$,
$E(G)$ and $\Delta(G)$ to denote the vertex set, edge set and maximum degree of a
graph $G$, respectively. For a vertex $v\in V(G)$, $N(v)$ is the set of vertices
adjacent to $v$, and $d(v)=|N(v)|$ is the {\em degree} of $v$. Similarly, we can
define $N(e)$, the set of edges adjacent to $e$. A vertex of degree $k$ (at least $k$,
or at most $k$) is called a {\em$k$-vertex} (a {\em $k^+$-vertex}, or a {\em
$k^-$-vertex}, respectively). For a vertex subset $S$ of $V(G)$, we use $G[S]$ to
denote the subgraph of $G$ that is induced by $S$. Let $n$, $m$ be two integers. A
complete bipartite graph with one part having $n$ vertices and the other part $m$
vertices is denoted by $K_{n,m}$. A graph is called {\em claw-free} if it has no
induced subgraph isomorphic to $K_{1,3}$.

\par An injective $k$-edge-coloring of a graph $G$ is a mapping $\phi$:
$E(G)\rightarrow\{1,2,...,k\}$, such that $\phi(e)\ne\phi(e')$ if edges $e$ and $e'$
are at distance two, or are in a triangle. The smallest integer $k$ such that $G$ has
an injective $k$-edge-coloring is called the {\em injective chromatic index} of $G$,
denoted by $\chi_i'(G)$. The concept of injective edge-coloring was proposed in 2015
by Cardose et al. \cite{DJJC-2019} to slove the Packet Radio Network problem and they
proved that it is NP-hard to compute the injective chromatic index for any graph.
Moreover, Ferdjallah et al.\cite{BSA-2022} showed that $\chi_i'(G)\le2(\Delta(G)-1)^2$
for any graph $G$ with $\Delta(G)\ge 3$; and $\chi_i'(G)\le 30$ for any planar graph
$G$. In particular, they proposed the following conjecture.

\begin{conjecture}
For every subcubic graph $G$, $\chi_i'(G)\le 6$.
\end{conjecture}

In 2022, Miao et al. \cite{ZYG-2022} posed the following conjecture.

\begin{conjecture}
For every simple graph $G$ with maximum degree $\Delta$, $\chi_i'(G)\le
\Delta(\Delta-1)$.
\end{conjecture}

\par Several authors have attacked this upper bound on the injective chromatic index
for graphs with small maximum degree. Towards Conjecture $1.1$,  Kostochka et al.
\cite{AAJ-2021} confirmed that $\chi_i'(G)\le 7$ for subcubic graphs and proved that
the upper bound $7$ can be improved to $6$ for subcubic planar graphs.

\par For graphs with maximum degree 4, we summarize the upper bounds of injective
chromatic index for graphs with maximum average degree restrictions.

\begin{theorem}

Let $G$ be a graph with $\Delta(G)=4$. We say the graph $G$ is a $(m,k)$-graph if
$\mathrm{mad}(G)<m$ and $\chi_i'(G)\le k$.\par

$(1)$ $G$ is a $(m,k)$-graph for $m=\frac{7}{3}$ and $k=5$ {\rm
\cite{XG-2025}};\par

$(2)$ $G$ is a $(m,k)$-graph for $(m,k)\in
\{(\frac{5}{2},6),(\frac{13}{5},7),(\frac{36}{13},8)\}$ {\rm \cite{JJ-2025}};\par

$(3)$ $G$ is a $(m,k)$-graph for
$(m,k)\in\{(\frac{14}{5},9),(3,10),(\frac{19}{6},11)\}$ {\rm \cite{ZYG-2022}};\par

$(4)$ $G$ is a $(m,k)$-graph for $m=\frac{33}{10}$ and $k=12$ {\rm
\cite{JX-2023}};\par

$(5)$ $G$ is a $(m,k)$-graph for $(m,k)\in\{(\frac{10}{3},13),(\frac{18}{5},14),(\frac{15}{4},15)\}$
{\rm \cite{YC-2018}}.\par

\end{theorem}

\par For claw-free graphs, Dong et al. \cite{DLL-2023} confirmed that the injective
chromatic index of any claw-free subcubic graph is less than or equal to 6 and the
upper bound 6 is tight in 2023. Cui and Han \cite{QZ-2024} proved that $\chi_i'(G)\le
5$ for every connected claw-free subcubic graph $G$ that is not isomorphic to $K_4$
and $\overline C_6$ in 2024.

\medskip

\par In this paper, we consider the injective chromatic index of claw-free graphs with
maximum degree at most 4.

\begin{theorem}
Let $G$ be a claw-free graph with $\Delta(G)\le4$. Then $\chi_i'(G)\le 13$.
\end{theorem}

\par Suppose that $G$ has a partial injective edge-coloring $\phi$ with the color set
$C$. For each edge $e'$ and $e$ in $G$, we say that edge $e'$ {\em sees} the edge $e$
if they are at distance two or are in a triangle. For $e=uv\in E(G)$, we denote the
set of the colors of the edges that see $e$ as $F_\phi(e)$ and denote the set of
available colors of $e$ as $S_\phi(e)$. Obviously, $S_\phi(e)=C-F_\phi(e)$ and
$|F_\phi(e)|\le 3(d(u)+d(v)-2)$. We simply write $S_\phi(e)$ as $S(e)$ if there is no
confusion. For a vertex $v\in V(G)$, we denote the set of the colors of the edges
incident with $v$ as $C_\phi(v)$.

For all figures in this paper, a vertex is represented by a solid point when all of
its incident edges are drawn; otherwise it is represented by a hollow point.
We will use the labels as shown in the figures.

\section{Proof of Theorem 1.2}
\par Assume that $G$ is a counterexample of Theorem $1.2$ such that $|V(G)|$ is as
small as possible. Recall that $\Delta(G)\leq 4$. Then $G$ is a connected claw-free
graph.

\begin{remark} Let $v\in V(G)$ and $uv\in E(G)$. Suppose that $u$ is not adjacent to
any other vertices in $N(v)\setminus\{u\}$. Since $G$ is claw-free, we have $xy\in
E(G)$ for any two vertices $x,y\in N(u)\setminus\{v\}$. So $vu$ sees at most $6$ edges
at the vertex $u$.
\end{remark}

\begin{lemma}
\par $\delta(G)=4$.
\end{lemma}
\begin{pf}
 \rm{Suppose to the contrary that $G$ contains a $3^-$-vertex $v$. Let $d(v)=k\le3$
 and $N(v)=\{v_1,v_2,\dots,v_k\}$. By the minimality of $G$, $G'=G-v$ has an injective
 13-edge-coloring $\phi$.

\noindent\textbf{Case$\;$1.}~$k=1$.

\par Since $|S(vv_1)|\ge 13-3(d(v_1)-1)\ge 4$, we can extend $\phi$ to $G$, a
contradiction.

\noindent\textbf{Case$\;$2.}~$k=2$.

\par First suppose that $v_1v_2\in E(G)$. Then $|S(vv_1)|\ge
13-3(d(v_1)-2)-(d(v_2)-1)\ge 4$ and $|S(vv_2)|\ge 13-3(d(v_2)-2)-(d(v_1)-1)\ge 4$, we
can extend $\phi$ to $G$, a contradiction.
\par Next suppose that $v_1v_2\notin E(G)$. Then $vv_1$ sees at most 6 edges at the
vertex $v_1$ by Remark 2.1. So $|S(vv_1)|\ge 13-(6+(d(v_2)-1))\ge4$. By symmetry,
$|S(vv_2)|\ge 4$. We can extend $\phi$ to $G$, a contradiction.

\noindent\textbf{Case$\;$3.}~$k=3$.

\par Set $q=|E(G[\{v_1,v_2,v_3\}])|$. Then $1\leq q\leq 3$ by $G$ is claw-free.

\noindent\textbf{Subcase$\;$3.1.}~$q=1$, say $v_1v_2\in E(G)$.

\par Then $v_2v_3\notin E(G)$ and $v_1v_3\notin E(G)$. Then $vv_3$ sees at most 6
edges at the vertex $v_3$  by Remark 2.1. Hence $|S(vv_3)|\ge
13-(6+(d(v_1)+d(v_2)-3))\ge2$. Next we can show that $|S(vv_1)|\ge 2$. In fact, if
$d(v_1)=3$ or $d(v_2)=3$, then $|S(vv_1)|\ge 13-3(d(v_1)-2)+(d(v_2)-1)+3)\ge 2$. Now
we can suppose that $d(v_1)=d(v_2)=4$. Then $xy\in E(G)$ by $G[\{v,x,y\}]$ is not
isomorphic to $K_{1,3}$, where $x,y\in N(v_1)\setminus \{v,v_2\}$. So $|S(vv_1)|\ge
13-(5+(d(v_2)-1)+(d(v_3)-1))\ge2$. Hence, we show that $|S(vv_1)|\ge 2$. By symmetry,
$|S(vv_2)|\ge 2$. Then $\phi$ can be extended to be an injective 13-edge-coloring of
$G$, a contradiction.

\noindent\textbf{Subcase$\;$3.2.}~$q=2$, say $v_1v_2\in E(G)$ and $v_2v_3\in E(G)$.
\par Then $|S(vv_1)|\ge 13-(3(d(v_1)-2)+(d(v_2)-1)+(d(v_3)-2))\ge 2$, and
$|S(vv_2)|\ge 13-(3+(d(v_1)-1)+(d(v_3)-1))\ge 4$. By symmetry,  $|S(vv_3)|\ge 2$. So
we can extend $\phi$ to $G$, a contradiction.

\noindent\textbf{Subcase$\;$3.3.}~$q=3$, say $v_1v_2\in E(G)$, $v_2v_3\in E(G)$ and
$v_1v_3\in E(G)$.
\par Then $|S(vv_i)|\ge 13-(3+5)=5$ for each $i\in \{1,2,3\}$, we can extend $\phi$ to
$G$, a contradiction.
\qed}
\end{pf}

\begin{lemma}
$G$ does not contain $K_4$ as a subgraph.
\end{lemma}

\begin{pf}
\rm{
\par Suppose that $G$ contains $K_4$ as a subgraph. Set $V(K_4)=\{v_1,v_2,v_3, v_4\}$.
Let $u_i$ be the neighbor of $v_i$ not in $V(K_4)$ for each $i\in\{1,2,3,4\}$.
\par Suppose that $u_1=u_2$. By the minimality of $G$, $G'=G-v_1$ has an injective
13-edge-coloring $\phi$. Since $|S(v_1u_1)|\ge
13-2(d(u_1)-2)-(d(v_2)+d(v_3)+d(v_4)-6)=1$, $|S(v_1v_2)|\ge
13-(d(u_1)-1)-(d(v_3)+d(v_4)-3)=5$, and $|S(v_1v_3)|\ge
13-(3+(d(v_4)+d(v_2)-3)+(d(u_1)-2))=3$. By symmetry, $|S(v_1v_4)|\ge 3$. So we can
extend $\phi$ to $G$, a contradiction.

\par So we may assume that any two of $u_1,u_2,u_3,u_4$ are not coincide. By the
minimality of $G$, $G'=G-\{v_1,v_2,v_3,v_4\}$ has an injective 13-edge-coloring
$\phi$. Then $v_iu_i$ sees at most 6 edges at the vertex $u_i$ by Remark 2.1 for each
$i\in \{1,2,3,4\}$. So $|S(u_iv_i)|\ge 13-6=7$. Since $|S(v_iv_j)|\ge
13-(d(v_i)+d(v_j)-2)=7$ for each pair $i,j\in \{1,2,3,4\}$, we can extend $\phi$ to
$G$ by coloring $v_1v_4,v_1v_2,v_2v_3,v_3v_4,v_1v_3,v_2v_4,v_4u_4,v_2u_2,v_3u_3$ and
$v_1u_1$ in order, a contradiction.
\qed
}
\end{pf}

\begin{lemma}
Any $4$-vertex in $G$ is incident with at most two 3-cycles.
\end{lemma}

\begin{pf}
\rm{

\par Suppose to the contrary that there exists $4$-vertex $v$ incident with three
3-cycles in $G$. Let $N(v)=\{v_1, v_2, v_3, v_4$\}. By Lemma 2.2 and $G$ is claw-free,
we may assume that $v_1v_2\in E(G)$, $v_2v_3\in E(G)$ and $v_1v_4\in E(G)$. Set
$N(v_1)=\{v,v_2,v_4,u_1\}$ and $N(v_2)=\{v,v_1,v_3,u_2\}$. By Lemma 2.2, $u_1\neq v_3$
and $u_2\neq v_4$. By the minimality of $G$, $G'=G-v$ has an injective
13-edge-coloring $\phi$. Since $G$ is claw-free, $u_1=u_2$, or $u_1v_4\in E(G)$ and
$u_2v_3\in E(G)$.

\par Suppose that $u_1=u_2$. First suppose that $v_3v_4\notin E(G)$. We erase the
color of $v_1v_4$. Then $vv_i$ sees at most 7 edges at the vertex $v$ for each
$i\in\{1,2\}$, and $vv_i$ sees at most 6 edges at the vertex $v$ for each
$i\in\{3,4\}$. Since $|S(vv_3)|\ge 13-(2\times3+6)=1$, $|S(vv_4)|\ge
13-(2\times3+6)=1$, $|S(v_1v_4)|\ge 13-(2\times3+4)=3$, $|S(vv_1)|\ge 13-(7+2)=4$ and
$|S(vv_2)|\ge 13-(7+2)=4$, we can extend $\phi$ to $G$, a contradiction. Next suppose
that $v_3v_4\in E(G)$. Then $vv_i$ sees at most 7 edges at the vertex $v$ for each
$i\in\{1,2,3,4\}$. Since $|S(vv_3)|\ge 13-(7+3)=3$, $|S(vv_4)|\ge 13-(7+3)=3$,
$|S(vv_1)|\ge 13-(7+2)=4$ and $|S(vv_2)|\ge 13-(7+2)=4$, we can extend $\phi$ to $G$,
a contradiction.

\par Suppose that $u_1v_4\in E(G)$ and $u_2v_3\in E(G)$. By the above case, we can
deduce that $v_3v_4\notin E(G)$. Then $vv_i$ sees at most 7 edges at the vertex $v$
for each $i\in\{3,4\}$, and $vv_i$ sees at most 8 edges at the vertex $v$ for each
$i\in\{1,2\}$. Since $|S(vv_3)|\ge 13-(7+3+2)=1$, $|S(vv_4)|\ge 13-(7+3+2)=1$,
$|S(vv_1)|\ge 13-(8+2)=3$ and $|S(vv_2)|\ge 13-(8+2)=3$, we can extend $\phi$ to $G$,
a contradiction.
\qed
}
\end{pf}

\par By Lemma 2.2 and 2.3, the following lemma holds trivially.

\begin{lemma}
\par Each $4$-vertex in $G$ is incident with exactly two edge-disjoint triangles.

\end{lemma}

\begin{figure}[htbp]
  \centering
  \includegraphics[width=1.7in]{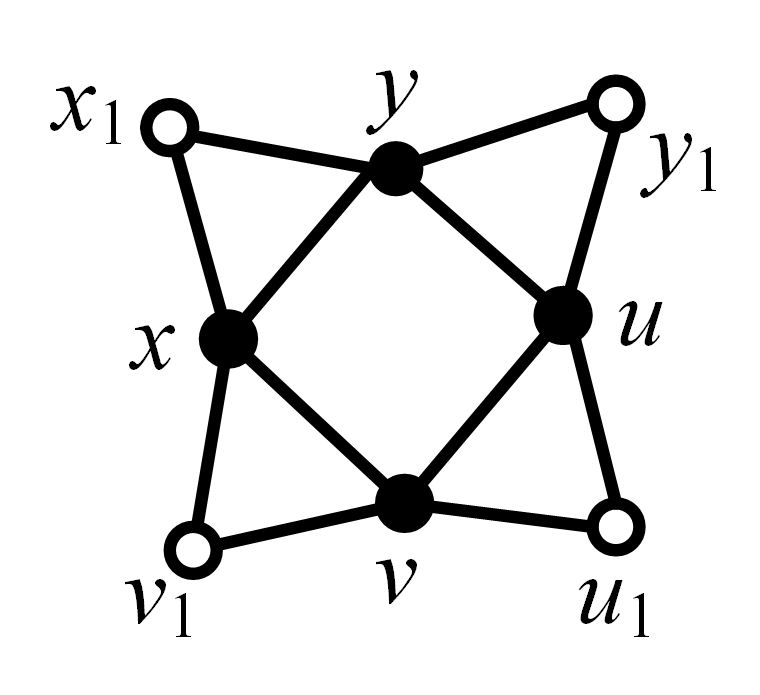}
  \caption{There exists a $4$-cycle $xyuvx$ in $G$.}
\end{figure}

\begin{lemma}
\par There is no $4$-cycles in $G$.

\end{lemma}
\begin{pf}
\rm{
\par Suppose to the contrary that there exists a $4$-cycle $xyuvx$. By Lemma 2.4, each
vertex in $\{x,y,u,v\}$ is incident with exactly two edge-disjoint triangles, as shown
in Figure 1.  Let $N(y_1)=\{y,u,y_1',y_1''\}$ and $N(u_1)=\{u,v,u_1',u_1''\}$. By Lemma 2.4, $y_1'y_1''\in E(G)$ and $u_1'u_1''\in E(G)$. By the
minimality of $G$, $G'=G-\{y,y_1,u,v,u_1\}$ has an injective 13-edge-coloring $\phi$.

\par
\textbf{Case 1}~$(N(y_1)\cap N(u_1))\setminus\{u\}\neq \emptyset$, say $y_1''=u_1'$.
\par We have $|S(yu)|\ge 13-(d(x_1)+d(x)-4)=9$. By
symmetry, $|S(uy_1)|\ge 9$, $|S(uv)|\ge 9$ and $|S(uu_1)|\ge 9$. Since $|S(yu)|+|S(uv)|\ge 18>13$, we have $|S(yu)\cap S(uv)|\ge 1$. We color $yu$ and $uv$ with a color $\alpha\in
S(yu)\cap S(uv)$, and denote this new coloring as $\phi'$. Now $|S_{\phi'}(uy_1)|\ge
9-1=8$ and $|S_{\phi'}(uu_1)|\ge 9-1=8$. Similarly, we can color $uy_1$ and $uu_1$
with a color $\beta \in S_{\phi'}(uy_1)\cap S_{\phi'}(uu_1)$ by
$|S_{\phi'}(uy_1)|+|S_{\phi'}(uu_1)|\ge 16>13$. Denote this new coloring as $\phi''$.
Then $|S_{\phi''}(yx_1)|\ge 13-((d(x)-2)+5+2)=4$, $|S_{\phi''}(yx)|\ge 13-((d(x_1)-1)+(d(v_1)-2)+2)=6$ and $|S_{\phi''}(yy_1)|\ge 13-((d(x_1)+d(x)-4)+(d(y_1')+d(y_1'')-4)+2)=3$. By symmetry, $|S_{\phi''}(y_1y_1')|\ge 4$, $|S_{\phi''}(y_1y_1'')|\ge 6$, $|S_{\phi''}(vv_1)|\ge 4$, $|S_{\phi''}(vx)|\ge 6$, $|S_{\phi''}(u_1u_1'')|\ge 4$, $|S_{\phi''}(u_1y_1'')|\ge 6$ and $|S_{\phi''}(vu_1)|\ge 3$. Hence we can extend $\phi''$ to $G$ by coloring $y_1y_1',y_1y_1'', u_1y_1'', u_1u_1'', yx_1, vv_1, yx, xv, yy_1$ and $vu_1$ in order, a contradiction.

\par
\textbf{Case 2}~$(N(y_1)\cap N(u_1))\setminus\{u\}=\emptyset$.

We have $|S(yu)|\ge 13-(d(x_1)+d(x)-4)=9$, $|S(uy_1)|\ge 13-(d(y_1')+d(y_1'')-3)=8$. By
symmetry, $|S(uv)|\ge 9$ and $|S(uu_1)|\ge 8$. Since $|S(yu)|+|S(uv)|\ge18>13$, we
have $|S(yu)\cap S(uv)|\ge 1$. We color $yu$ and $uv$ with a color $\alpha\in
S(yu)\cap S(uv)$, and denote this new coloring as $\phi'$. Now $|S_{\phi'}(uy_1)|\ge
8-1=7$ and $|S_{\phi'}(uu_1)|\ge 8-1=7$. Similarly, we can color $uy_1$ and $uu_1$
with a color $\beta \in S_{\phi'}(uy_1)\cap S_{\phi'}(uu_1)$ by
$|S_{\phi'}(uy_1)|+|S_{\phi'}(uu_1)|\ge14>13$. Denote this new coloring as $\phi''$.
Then $|S_{\phi''}(y_1y_1')|\ge 13-((d(y_1'')-1)+5+2)=3$, $|S_{\phi''}(yx_1)|\ge
13-((d(x)-2)+5+2)=4$, $|S_{\phi''}(yx)|\ge 13-((d(x_1)-1)+(d(v_1)-2)+2)=6$ and
$|S_{\phi''}(yy_1)|\ge 13-((d(x_1)-1)+(d(x)-3)+(d(y_1')+d(y_1'')-3)+2)=2$. By
symmetry, $|S_{\phi''}(y_1y_1'')|\ge 3$, $|S_{\phi''}(u_1u_1')|\ge 3$,
$|S_{\phi''}(u_1u_1'')|\ge 3$, $|S_{\phi''}(vv_1)|\ge 4$, $|S_{\phi''}(xv)|\ge 6$ and
$|S_{\phi''}(vu_1)|\ge 2$. Hence we can extend $\phi''$ to $G$ by coloring $y_1y_1',
y_1y_1'', u_1u_1', u_1u_1'', yx_1, vv_1, yx, xv, yy_1$ and $vu_1$ in order, a
contradiction.\qed

}
\end{pf}

\rm{
\par Now we are ready to show Theorem 1.2. Let $v\in V(G)$ with
$N(v)=\{u_1,u_2,u_3,u_4\}$. By Lemma 2.4, we may assume that $u_1u_2\in E(G)$ and
$u_3u_4\in E(G)$. Let $N(u_1)=\{x_1,y_1,v,u_2\}$, $N(u_2)=\{x_2,y_2,v,u_1\}$,
$N(u_3)=\{x_3,y_3,v,u_4\}$, and $N(u_4)=\{x_4,y_4,v,u_3\}$. By Lemma 2.4, $x_iy_i\in
E(G)$, $u_i\neq x_j$ and $u_i\neq y_j$ for each pair $i,j\in \{1,2,3,4\}$. By Lemma
2.5, $x_i\neq x_j$, $x_i\neq y_i$ for each $i,j\in\{1,2,3,4\}$ and $i\neq j$,  as shown in Figure 2.
By the minimality of $G$, $G'=G-\{v,u_1,u_2,u_3,u_4\}$ has an injective
13-edge-coloring $\phi$ with the color set $C$.

\begin{figure}[htbp]
  \centering
  \includegraphics[width=2.5in]{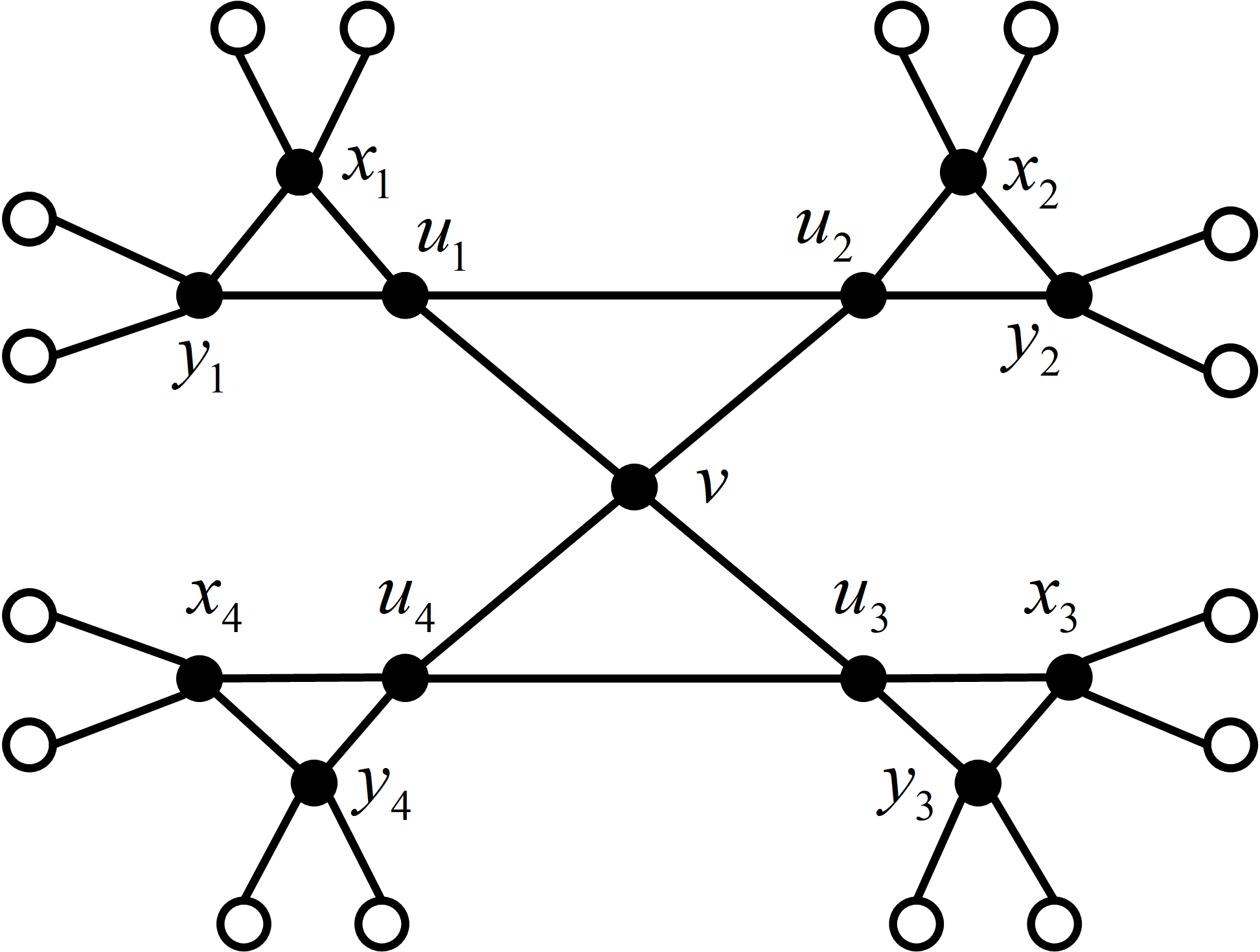}
  \caption{The configuration used in the proof of Theorem 1.2.}
\end{figure}

\par Since $G$ is claw-free and $\Delta(G)\leq 4$, we have $|S(vu_i)|\ge 13-5=8$,
$|S(u_ix_i)|\ge 13-8=5$, $|S(u_iy_i)|\ge 13-8=5$ for each $i\in \{1,2,3,4\}$. Since
$|S(vu_1)|+|S(vu_3)|=16>13$, we have $|S(vu_1)\cap S(vu_3)|\ge 1$. We color $vu_1$ and
$vu_3$ with a color $\alpha\in S(vu_1)\cap S(vu_3)$, and denote this new coloring as
$\phi'$. Now $|S_{\phi'}(u_ix_i)|\ge 5-1=4$ and $|S_{\phi'}(u_iy_i)|\ge 5-1=4$ for
each $i\in \{1,2,3,4\}$. We can color
$u_1x_1,u_1y_1,u_2x_2,u_2y_2,u_3x_3,u_3y_3,u_4x_4$ and $u_4y_4$ with
$b_1,b_2,b_3,b_4,b_5,b_6,b_7,b_8$ in order, and denote the obtained new coloring as
$\phi''$. Then $|S_{\phi''}(vu_2)|\ge 8-1-6=1$ and $|S_{\phi''}(v u_4)|\ge 8-1-6=1$.
We color $vu_2$ with $\beta_1$, and color $vu_4$ with $\beta_2$. Denote the
obtained new coloring as $\phi^*$.

\par Suppose that $\beta_1=\beta_2$. Since $|S_{\phi^*}(u_1u_2)|\ge 13-(5+5+2)=1$ and
$|S_{\phi^*}(u_3u_4)|\ge 13-(5+5+2)=1$, we can extend $\phi^*$ to $G$, a
contradiction. So we may assume that $\beta_1\neq \beta_2$, say $\alpha=1,\beta_1=2$
and $\beta_2=3$. If $|S_{\phi^*}(u_1u_2)|\ge 1$ and $|S_{\phi^*}(u_3u_4)|\ge 1$, then
we can extend $\phi^*$ to $G$, a contradiction. By symmetry, we may assume that
$|S_{\phi^*}(u_1u_2)|=0$. That is $(C_{\phi^*}(x_1)\cup C_{\phi^*}(y_1))\setminus
\{b_1,b_2\}=\{4,5,6,7,8\}$ and $(C_{\phi^*}(x_2)\cup C_{\phi^*}(y_2))\setminus
\{b_3,b_4\}=\{9,10,11,12,13\}$. If we can recolor $vu_2$ with 3, then turn to the case
that $\beta_1=\beta_2$, a contradiction. Hence $3\in F_{\phi*}(vu_2)$, say
$3\in\{b_1,b_2,b_5,b_6,b_7,b_8\}$. We can deduce that $3\in \{b_7,b_8\}$, say $b_7=3$,
by $\phi^*$ is the partial injective edge-coloring of $G$.

\begin{itemize}
\item
If there exists a color $\gamma\in\{4,5,6,7,8\}$ such that $\gamma\notin
F_{\phi*}(vu_2)$, say $\gamma=4$, then we recolor $vu_2$ with $\gamma$ and color
$u_1u_2$ with 2. Denote this new coloring as $\phi_1$. Now if
$|S_{\phi_1}(u_3u_4)|\ge1$, then we can extend $\phi_1$ to $G$, a contradiction. Hence
$F_{\phi_1}(u_3u_4)=C$. Let $(C_{\phi_1}(x_3)\cup C_{\phi_1}(y_3))\setminus
\{b_5,b_6\}=\{c_1,c_2,\ldots,c_5\}$ and $(C_{\phi_1}(x_4)\cup
C_{\phi_1}(y_4))\setminus \{b_7,b_8\}=\{d_1,d_2,\ldots,d_5\}$. Then
$\{c_1,\dots,c_5,d_1,\dots,d_5\}=\{2,5,6,\dots,13\}$. If we can recolor $vu_4$ with
$4$, then we can color $u_3u_4$ with 3 to obtain an injective 13-edge-coloring of
$G$, a contradiction. Hence $4\in F_{\phi_1}(vu_4)$, say $4\in
\{b_1,b_2,\ldots,b_6\}$. By $\phi_1$ is the partial injective edge-coloring of $G$,
we have $4\in \{b_3,b_4\}$, say $b_3=4$. Note that $3\notin S_{\phi_1}(vu_3)$. If we
can recolor $vu_3$ with a color $\gamma\in S_{\phi_1}(vu_3)\setminus\{1\}$, then we
can recolor or color $vu_4,u_3u_4$ with 1, 3, respectively. The obtained coloring is
the injective 13-edge-coloring of $G$, a contradiction. Hence $F_{\phi_1}(vu_3)\cup
\{1\}=C$. That is $\{b_1,b_2,b_4,b_8,2,c_1,c_2,\dots,c_5\}=\{2,5,6,\dots,13\}$.
Recall that $\{c_1,c_2,\dots,c_5,d_1,d_2,\dots,d_5\}=\{2,5,6,\dots,13\}$. We can
deduce that $\{b_1,b_2,b_4,b_8\}\subseteq\{d_1,d_2,\dots,d_5\}$. Now
$|F_{\phi_1}(vu_4)|\le 10$, we can recolor $vu_4$ with a color $\eta\in
S_{\phi_1}(vu_4)\setminus\{3\}$ and color $u_3u_4$ with 3 to obtain the injective
13-edge-coloring of $G$, a contradiction.

\item
 If $\{4,5,6,7,8\}\subseteq F_{\phi*}(vu_2)$, then
 $\{b_1,b_2,b_5,b_6,b_8\}=\{4,5,6,7,8\}$.
Now we erase the color of $vu_1$, and recolor or color $vu_2,u_1u_2$ with 1, 2,
respectively. We denote this new coloring as $\phi_2$. Since
$|S_{\phi_2}(u_3u_4)|\ge 13-12=1$ and $|S_{\phi_2}(vu_1)|\ge 13-10=3$, we can extend
$\phi_2$ to $G$, a contradiction.
\end{itemize}
}


\begin{thebibliography}{99}
\bibitem{YC-2018} Y. Bu, C. Qi, Injective edge coloring of sparse graphs, Discrete
    Math. Algorithms Appl. 10 (2018) 1850022, 16 pp.
\bibitem{DJJC-2019} D. Cardoso, J. Cerdeira, J. Cruz, C. Dominic, Injective edge
    coloring of graphs, Filomat 33 (2019) 6411-6423.
\bibitem{QZ-2024} Q. Cui, Z. Han, Injective edge-coloring of claw-free subcubic
    graphs, J. Comb. Optim. 47 (2024) 87, 32 pp.
\bibitem{DLL-2023} X. Dong, Y. Lin, W. Lin, The injective chromatic index of a
    claw-free subcubic graph is at most 6, J. Math. Res. Appl. 43 (2023) 409-416.
\bibitem{BSA-2022} B. Ferdjallah, S. Kerdjoudj, A. Raspaud, Injective edge-coloring of
    subcubic graphs, Discrete Math. Algorithms Appl. 14 (2022) 2250040, 22 pp.
\bibitem{JJ-2025} J. Fu, J. Lv, On injective edge-coloring of graphs with maximum
    degree 4, Discrete Appl. Math. 360 (2025) 119-130.
\bibitem{XG-2025} X. Hu, G. Zhang, Injective edge chromatic index of sparse graphs, Discrete Appl. Math. 370 (2025) 50-56.
\bibitem{AAJ-2021} A. Kostochka, A. Raspaud, J. Xu, Injective edge-coloring of graphs
    with given maximum degree, European J. Combin. 96 (2021) 103355, 12 pp.
\bibitem{JX-2023} J. Lu, X. Pan, Injective edge coloring of some sparse graphs, J.
    Appl. Math. Comput. 69 (2023) 3421-3431.
\bibitem{ZYG-2022} Z. Miao, Y. Song, G. Yu, Note on injective edge-coloring of graphs,
    Discrete Appl. Math. 310 (2022) 65-74.
\end{thebibliography}
\end{document}